\documentclass[11pt, epsf]{amsart}
\usepackage{epsfig}
\usepackage{eucal}
\usepackage{mathrsfs}
\usepackage{amssymb}
\usepackage{latexsym}
\usepackage{amsfonts}
\usepackage[latin1]{inputenc}
\usepackage{color}

\newenvironment{namelist}[1]{%
\begin{list}{}
 {
   
   \settowidth{\labelwidth}{#1}
   \setlength{\leftmargin}{1.1\labelwidth}
  }
 }{%
\end{list}}
\newcommand{\bi}{\begin{namelist}}
\newcommand{\ei}{\end{namelist}}

\newcommand{\vt}{\vartriangle}
\newcommand{\dis}{\displaystyle}
\newcommand{\eps}{\varepsilon}

\newtheorem{Theo}{Theorem}

\newtheorem{Lem}{Lemma}
\newtheorem{Rem}{Remark}

\title[Strong convergence of wave operators for Dirac operators]{Strong convergence of wave operators for a family of Dirac operators}
\author{Hasan Almanasreh}
\thanks{1- Mathematics Department, Hebron University,
        P.O. Box 40, Hebron, West Bank, Palestine}
\thanks{2- Department of Math. Sciences, University of Gothenburg,
        SE-412 96 Gothenburg, Sweden}
\thanks{E-mail: hasanm@hebron.edu}

\begin{document}

\maketitle

\begin{abstract} We consider a family of Dirac operators with potentials varying with respect to a parameter $h$. The set of potentials considered has different power-like decay that is independent of $h$. The proofs of existence and completeness of the wave operators are similar to that given in \cite{GAT}. We are mainly interested in the asymptotic behavior of the wave operators as $h\to\infty$.\\
\end{abstract}
\section{Introduction}
In quantum mechanics, it is important to compare a given interacted operator with a simpler (free) operator for which many spectral features are known. Scattering theory is part of perturbation theory that concerns a comparative study of the absolutely continuous spectrum of operators. That is, for two self-adjoint operators $T$ and $T_0$ that are close to each other in an appropriate sense, scattering theory is mainly the study of existence and completeness of $\dis s\!-\!\!\!\!\lim_{t\to\pm\infty} \dis e^{iTt}\mathcal{J}e^{-iT_0t}P_0^{(ac)}$ where $\mathcal{J}$ is some bounded operator (identification), $P_0^{(ac)}$ is the orthogonal projection onto the absolutely continuous spectrum of $T_0$, and $s$ refers to the strong convergence sense. Another important issue is the case studied in the present work, where the operator $T$ is $h$-dependent, and $h$ is a parameter allowed to grow to infinity. For such operators, in the addition to the existence study of the time limit, a parallel question also emerges whether or not the limit $\dis s\!-\!\!\!\lim_{h\to\infty}\dis s\!-\!\!\!\!\lim_{t\to\pm\infty} \dis e^{iT_ht}\mathcal{J}_he^{-iT_0t}P_0^{(ac)}$ exists, where now $\mathcal{J}_h$ is an $h$-dependent bounded identification.

Scattering theory for the Dirac operator with potentials decaying faster than the Coulomb potential (short-range potentials) goes straightforward. In this case, the proofs of existence and completeness of the wave operator (WO) $W(H,H_0)=\dis s\!\!-\!\!\!\!\!\lim_{t\to\pm\infty}e^{iHt}e^{-iH_0t}$, where $H_0$ and $H$ are respectively the free and the interacted Dirac operators, are similar to that of the Schr\"odinger operator. In the Coulomb interaction case, the modified WO $W_{\pm}=W_\pm(H,H_0,\mathcal{J})=\dis s\!\!-\!\!\!\!\!\lim_{t\to\pm\infty}e^{iHt}\mathcal{J}e^{-iH_0t}$ has been constructed in \cite{DOL,DOLV}, where $\mathcal{J}$ is a bounded identification. For potentials decaying as the Coulomb potential or slower (long-range potentials), the existence and completeness of the modified WO $W_\pm$ have been studied in \cite{GAT,MUT,MUTS,THAE,YAM}. The study of the asymptotic behavior of the WO $W_\pm$ with respect to the speed of light ($c$), as $c\to\infty$, has been studied for short-range potential in \cite{YAJ} and for long-range potential in \cite{YAM}.

Consider the family of Dirac operators $H_h=H_0+V_h$ defined on the same Hilbert space and with the same domain as of $H_0$, where $H_0$ is the free Dirac operator defined on $L^2(\mathbb {R}^3,\mathbb {C}^4)$ with domain $H^1(\mathbb {R}^3,\mathbb {C}^4)$, and $V_h$ is a bounded interaction to $H_0$. Under suitable power-like decay assumption on $V_h$, the WO $W_{\pm,h}=W_\pm(H_h,H_0;\mathcal{J}_{\pm,h})=\dis s\!-\!\!\!\!\lim_{t\to\pm\infty}e^{iH_ht}\mathcal{J}_{\pm,h}e^{-iH_0t}$ exists and is complete \cite{GAT}, where $\mathcal{J}_{\pm,h}$ is a bounded identification. In other words, if for all $h>0$, $|V_h|\leq\langle x\rangle^{-\rho}$, where $\langle x\rangle=(1+|x|^2)^{1/2}$, then the WO $W_{\pm,h}$ exists and is complete where $\mathcal{J}_{\pm,h}$ being just the identity operator for $\rho>1$ (short-range). For $0<\rho\leq1$ (long-range), the identification $\mathcal{J}_{\pm,h}$ is a pseudo-differential operator (PSDO) defined as
\begin{equation*}
(\mathcal{J}_{\pm,h}g)(x)=(2\pi)^{-3/2}\dis\int_{\mathbb{R}^3}\dis e^{i x\cdot\zeta + i\Phi_{\pm,h}(x,\zeta)} \mathcal{P}_{\pm,h}(x,\zeta)\mathscr{C}_{\pm}(x,\zeta)\psi(|\zeta|^2)\hat{g}(\zeta) \,d\zeta,
\end{equation*}
where $\Phi_{\pm,h}$ is a phase function, $\mathcal{P}_{\pm,h}$ is an amplitude function, $\mathscr{C}_{\pm}$ is a cut-off function, and $\psi$ is a smooth function introduced to localize $\mathcal{J}_{\pm,h}$ in compact intervals of the continuous spectrum.

The goal of the present work is to study the asymptotic behavior of the WO $W_{\pm,h}$ and its adjoint $W_{\pm,h}^*=W_\pm(H_0,H_h;\mathcal{J}_{\pm,h}^*)$ as $h\to\infty$. By the existence of $W_{\pm,h}$, the convergence of $H_h$ in the strong resolvent sense (SRS), and the strong convergence of the identification $\mathcal{J}_{\pm,h}$, we prove that the two limits $\dis s\!-\!\!\!\lim_{h\to\infty}$ and $\dis s\!-\!\!\!\!\lim_{t\to\pm\infty}$ are interchangeable. Hence, if the Dirac operator $H_h$ converges in the SRS to $H_\infty$ and $\mathcal{J}_{\pm,h}$ converges strongly to $\mathcal{J}_{\pm,\infty}$, then we have
\begin{equation*}
s\!-\!\!\!\dis\lim_{h\to\infty}W_{\pm,h}= W_{\pm}(H_\infty,H_0;\mathcal{J}_{\pm,\infty}).
\end{equation*}
By the strong convergence of $\mathcal{J}_{\pm,h}$ to $\mathcal{J}_{ \pm,\infty}$, the identification $\mathcal{J}_{\pm,h}^*$ also converges strongly to $\mathcal{J}_{ \pm,\infty}^*$ (this is not true in general, but it is valid for the type of identifications we consider in this work), hence we also have
\begin{equation*}
s\!-\!\!\!\dis\lim_{h\to\infty}W_{\pm,h}^*=W_{\pm}(H_0,H_\infty;\mathcal{J}_{\pm,\infty}^*).
\end{equation*}

\section{Preliminaries}
By $\mathbf{R}(A)$, $\mathbf{D}(A)$, and $\mathbf{N}(A)$ we refer respectively to the range, domain, and null space of a given operator $A$, also we denote by $X$ and $Y$ the Hilbert spaces $H^1(\mathbb {R}^3,\mathbb {C}^4)$ and $L^2(\mathbb {R}^3,\mathbb {C}^4)$ respectively.
\subsection{The Dirac operator with an $h$-dependent potential}

The free Dirac evolution equation is given by
\begin{equation}\label{1}
i\hslash\frac{\partial}{\partial t}u(x,t) =H_0u(x,t),\quad\quad u(x,0)=u^0(x),
\end{equation}
where $H_0:  X\longrightarrow Y$ is the free Dirac operator defined as
\begin{equation}\label{2}
H_0 = \hslash cD_\alpha + mc^2 \beta.
\end{equation}
The constant $c$ is the speed of light, $\hslash$ is the Planck constant divided by $2\pi$, $D_\alpha=\alpha\cdot D$ where $D=(D_1,D_2,D_3)$ and $D_j=-i\frac{\partial}{\partial x_j}$ for $j=1,2,3$, and $m$ is the particle rest mass. The notations $\alpha  = (\alpha_1,\alpha_2,\alpha_3)$ and $\beta$ are the $4\!\times\!4$ Dirac matrices given by
$$
\alpha_j = \left(
\begin{array}{cc}
0 & \sigma_j \\
\sigma_j & 0
\end{array}
\right)\;\;\text{and}\;\;
\beta = \left(
\begin{array}{cc}
I & 0 \\
0 & -I
\end{array}
\right).
$$
Here, $I$ and $0$ are the $2\!\times\!2$ unity and zero matrices respectively, and $\sigma_j$'s are the $2\!\times\!2$ Pauli matrices
$$
\sigma_1 = \left(
\begin{array}{cc}
0 & 1 \\
1 & 0
\end{array}
\right),\;\;
\sigma_2 = \left(
\begin{array}{cc}
0 & -i \\
i & 0
\end{array}
\right)
\;\;\text{and}\;\;
\sigma_3 = \left(
\begin{array}{cc}
1 & 0 \\
0 & -1
\end{array}
\right).
$$
Separating the variables $x$ and $t$ in (\ref{1}) yields the free Dirac eigenvalue problem
\begin{equation}\label{3}
H_0u(x)=\lambda u(x),
\end{equation}
where $u(x)$ is the spatial part of the wave function $u(x,t)$, and $\lambda$ is the total energy of the particle. The free operator $H_0$ is essentially self-adjoint on $C^\infty_0(\mathbb{R}^3\backslash\left\{0 \right\},\mathbb{C}^4)$ and self-adjoint on $X$, its spectrum, $\sigma(H_0)$, is purely absolutely continuous and is given by
$$
\sigma(H_0) = (-\infty,-mc^2]\cup[mc^2,+\infty).
$$

Let $\mathscr{F}$ be the Fourier transform operator
\begin{equation}\label{4}
(\mathscr{F}u)(\zeta)=(2\pi)^{-3/2}\displaystyle \int_{\mathbb{R}^3} e^{-ix\cdot\zeta}u(x)dx=:\hat{u}(\zeta),
\end{equation}
then $\mathscr{F}H_0\mathscr{F}^*$ is the multiplication operator given by the matrix
\begin{equation}\label{5}
h_0(\zeta)=\zeta_\alpha + mc^2\beta,
\end{equation}
known as the symbol of $H_0$, where $\zeta_\alpha=\alpha\cdot\zeta=\displaystyle\sum_{k=1}^3\alpha_k\zeta_k$. The symbol $h_0(\zeta)$ can be written as
\begin{equation}\label{6}
h_0(\zeta)=\eta(\zeta)p_{+,0}(\zeta)-\eta(\zeta)p_{-,0}(\zeta),
\end{equation}
where $p_{\pm,0}(\zeta)$ are the orthogonal projections onto the eigenspaces of $h_0(\zeta)$ and are given by
\begin{equation}
p_{\pm,0}(\zeta)=\frac{1}{2}\big(I\pm\eta^{-1}(\zeta)(\zeta_\alpha+mc^2\beta)\big),
\end{equation}
and $\pm\eta(\zeta)=\pm\sqrt{|\zeta|^2 +m^2c^4}$ are the corresponding eigenvalues.

We consider an $h$-dependent potential, $V_h(x)$, added to the free Dirac operator and define
\begin{equation}\label{7}
H_h=H_0+V_h.
\end{equation}
The potential $V_h$ is assumed to be real and say bounded, thus, for all $h>0$, $H_h$ and $H_0$ have the same domain $X$ and $H_h$ is self-adjoint on $X$. For simplicity, we assume $\hslash=c=1$. The corresponding evolution equation reads
\begin{equation}\label{8}
\left\{ \begin{array}{l}
i\frac{\partial}{\partial t}u_h(x,t) =H_h u_h(x,t),\\
u_h(x,0)=u_h^0(x).
\end{array} \right.
\end{equation}
By the Stone theorem, there exists a unique solution to (\ref{8}) given by
\begin{equation}\label{9}
u_h(x,t)=U_h(t)u_h^0(x)\, , \; u_h^0\in X,
\end{equation}
where the strongly continuous unitary operator $U_h(t)=e^{-i H_ht}$ is generated by the operator $-i H_h$, see e.g. \cite{KAT76,THA}.

The potential $V_h$ is assumed to fulfill the following condition for all multi-index $\alpha$
\begin{equation}\label{10}
|\partial^\alpha V_{h}(x)|\leq C\langle x\rangle^{-\rho-|\alpha|},\quad \text{for all }h>0, \text{ and } \rho\in(0,1].
\end{equation}
The constant $C$ is independent of $x$ and $h$, and again $\langle x\rangle=(1+|x|^2)^{1/2}$. This condition simply means that $V_h$ is of long-range type for all $h>0$.

\subsection{Basic setting of scattering theory}
Given self-adjoint operators $T_0$ and $T$ in Hilbert spaces $\mathscr{H}_0$ and $\mathscr{H}$ respectively. Let $P_0^{(ac)}$ and $P^{(ac)}$ be the orthogonal projections onto the absolutely continuous subspaces, $\mathscr{H}_0^{(ac)}$ and $\mathscr{H}^{(ac)}$, of $T_0$ and $T$ respectively. The WO for $T$ and $T_0$, with a bounded identification $\mathcal{J}:\mathscr{H}_0\to \mathscr{H}$, denoted by $W_\pm(T,T_0;\mathcal{J})$, is defined as
\begin{equation}\label{15}
W_\pm(T,T_0;\mathcal{J})=\displaystyle s\!-\!\!\!\!\lim_{t\to\pm\infty}U(-t)\mathcal{J} U_0(t)P_0^{(ac)},
\end{equation}
provided that the corresponding strong limits exist, where again the letter $s$ refers to the strong convergence sense, $U(t)=\displaystyle e^{-iTt}$ and $U_0(t)=\displaystyle e^{-iT_0t}$. If $\mathscr{H}=\mathscr{H}_0$ and $\mathcal{J}$ is the identity operator, then the WO is denoted by $W_\pm(T,T_0)$. Also if $T_0$ has only absolutely continuous spectrum, then $P_0^{(ac)}$ is superfluous.

If the WO exists, then it is bounded. Since the operator $U(-t)U_0(t)$ is unitary, the WO $W_\pm(T,T_0)$ is isometric. In the case that $\mathcal{J}$ is not the identity operator, the WO $W_\pm(T,T_0;\mathcal{J})$ is isometric if $\displaystyle \lim_{t\to\pm\infty}\|\mathcal{J}U_0(t)u_0\|_\mathscr{H}$ $=\|u_0\|_{\mathscr{H}_0}$ for any $u_0\in \mathscr{H}_0^{(ac)}$. The WO admits the chain rule, i.e., if $W_\pm(T,T_1;\mathcal{J}_1)$ and $W_\pm(T_1,T_0;\mathcal{J}_0)$ exist, then the WO $W_\pm(T,T_0;\mathcal{J}_{10})=W_\pm(T,T_1;\mathcal{J}_1)W_\pm(T_1,T_0;\mathcal{J}_0)$ also exists, where $\mathcal{J}_{10}=\mathcal{J}_1\mathcal{J}_0$. The WO possesses the intertwining property, that is
\begin{equation}\label{16}
\phi(T)W_\pm(T,T_0;\mathcal{J})=W_\pm(T,T_0;\mathcal{J})\phi(T_0),
\end{equation}
for any bounded Borel function $\phi$. Also for any Borel set $\vt\subset \mathbb{R}$
\begin{equation}\label{17}
E(\vt)W_\pm(T,T_0;\mathcal{J})=W_\pm(T,T_0;\mathcal{J})E_0(\vt),
\end{equation}
where $E$ and $E_0$ are the spectral families of $T$ and $T_0$ respectively. The following remark is about the equivalence between WOs with different identifications.
\begin{Rem}
\emph{
Assume that the WO $W_\pm(T,T_0;\mathcal{J}_1)$ exists with an identification $\mathcal{J}_1$, and let $\mathcal{J}_2$ be another identification such that the difference $\mathcal{J}_1-\mathcal{J}_2$ is compact, then the WO $W_\pm(T,T_0;\mathcal{J}_2)$ exists and $W_\pm(T,T_0;\mathcal{J}_1)=W_\pm(T,T_0;\mathcal{J}_2)$. Moreover, the condition that $\mathcal{J}_1-\mathcal{J}_2$ is compact can be replaced by $s\!-\!\!\!\!\dis\lim_{t\to\pm\infty}(\mathcal{J}_1-\mathcal{J}_2)U_0(t)P_0^{(ac)}=0$.
}
\end{Rem}
The WO $W_\pm$ is said to be complete if $\mathbf{R}(W_\pm)=\mathscr{H}^{(ac)}$, and if $W_\pm(T,T_0;\mathcal{J})$ is complete, then the absolutely continuous part of $T_0$ is unitary equivalent to that of $T$. We refer to \cite{YAF92} for the completeness criteria. For comprehensive materials on scattering theory we refer to \cite{KAT67,KAT76,KAT71,REEvol3,WEI80,YAF92, YAF00}.
\begin{Rem}
\emph{
Given a self-adjoint operator $T$ in a Hilbert space $\mathscr{H}$. A $T$-bounded operator, $A:\mathscr{H}\to\mathfrak{H}$, where $\mathfrak{H}$ is an auxiliary Hilbert space, is called $T$-smooth if one of the following properties is fulfilled
$$
\dis\sup_{\|v\|_\mathscr{H}=1,v\in\mathbf{D}(T)}\dis\int_{-\infty}^\infty\|A \dis e^{-iTt}v\|_\mathfrak{H}^2\,dt<\infty.
$$
$$
\dis\sup_{\eps>0,\mu\in\mathbb{R}}\|A R_T(\mu\pm i\eps)\|^2<\infty,
$$
where $ R_T(\mu\pm i\eps)$ is the resolvent operator of $T$.
}
\end{Rem}
\subsection{The wave operator for a family of Dirac operators}
Suppose that, as $h\to\infty$, the Dirac operator $H_h$ converges to $H_\infty$ in the SRS, and $\mathcal{J}_{\pm,h}$ converges strongly to $\mathcal{J}_{\pm,\infty}$. Then the WOs $W_\pm^\dag$ and $W_\pm^{\dag,*}$ exist, and
\begin{equation*}
W_\pm^\dag=s\!-\!\!\!\dis\lim_{h\to\infty}W_{\pm}(H_h,H_0;\mathcal{J}_{\pm,h}) =W_{\pm}(H_\infty,H_0;\mathcal{J}_{\pm,\infty}).
\end{equation*}
\begin{equation*}
W_\pm^{\dag,*}=s\!-\!\!\!\dis\lim_{h\to\infty}W_{\pm}(H_0,H_h;\mathcal{J}_{\pm,h}^*) =W_{\pm}(H_0,H_\infty;\mathcal{J}_{\pm,\infty}^*).
\end{equation*}

The proof follows Theorem 2.1 in \cite{BRU}, and is divided into several steps. Firstly, see \cite{GAT},
\begin{equation}\label{36}
H_h\mathcal{J}_{\pm,h}-\mathcal{J}_{\pm,h}H_0=\dis\sum_{j=1}^3T_j^*B_{1,h}T_j+ \langle x\rangle^{-(1+\rho)/2}B_{2,h}\langle x\rangle^{-(1+\rho)/2},
\end{equation}
where $T_j=\langle x\rangle^{-1/2}\nabla_j$, $(\nabla_jv)(x)=\partial_jv(x)-|x|^2x_j \dis\sum_{k=1}^3x_k(\partial_kv)(x)$, $B_{1,h}$ and $B_{2,h}$ are bounded operators. Note that for all $h>0$, $\langle x\rangle^{-(1+\rho)/2}$ for all $\rho>0$ and $T_j$ for j=1,2,3, are $H_0$-smooth and $H_h$-smooth on any compact set $\vt\subset\!(-\infty,-m)\cup (m,\infty)$ such that $\vt\!\!\cap\sigma_p(H_h)=\emptyset$, where $\sigma_p$ is the point spectrum. The $H_0$-smoothness and $H_h$-smoothness of $\langle x\rangle^{-(1+\rho)/2}$ and of $T_j$ are known respectively as the limiting absorption principle and the radiation estimate. It is a fact that operator-smoothness is invariant under the multiplication by a bounded operator from left (or a bounded operator from right provided it commutes with the spectral family of the given operator), thus we can rewrite (\ref{36}) as
\begin{equation}\label{37}
H_h\mathcal{J}_{\pm,h}-\mathcal{J}_{\pm,h}H_0=\dis\sum_{i=1}^2K_{2,i,h}^*K_{1,i,h}.
\end{equation}
The operators $K_{2,1,h}$ and $K_{2,2,h}$ are $H_h$-smooth for all $h>0$, and $K_{1,1,h}$ and $K_{1,2,h}$ are $H_0$-smooth for all $h>0$. Without loss of generality, we may assume that, for all $h>0$,
\begin{equation}\label{38}
H_h\mathcal{J}_{\pm,h}-\mathcal{J}_{\pm,h}H_0=G_h^*G_{0,h}
\end{equation}
for an $H_h$-smooth and $H_0$-smooth operators $G_h$ and $G_{0,h}$ respectively.
\begin{Lem}
\emph{
For all $h>0$ and for all $u_0\in X$ the function
\begin{equation}\label{39}
\mathscr{K}^{(1)}_{u_0,h}(t)=\|(H_h\phi(H_h)\mathcal{J}_{\pm,h}\phi(H_0)-\phi(H_h) \mathcal{J}_{\pm,h}H_0\phi(H_0))U_0(t)u_0\|_Y
\end{equation}
belongs to $L^1((-\infty,\infty);dt)$ for some continuous function $\phi:\mathbb{R} \to\mathbb{R}$ such that $x\phi(x)$ is bounded on $\mathbb{R}$.
}
\end{Lem}
\hspace{-4mm}\underline{\emph{Proof}}. Let $\phi(x)=(x-z)^{-1}$, $z\in \operatorname{Res}(H_h)\cap \operatorname{Res}(H_0)$ where $\operatorname{Res}$ denotes the resolvent set. Therefore, and since
\begin{equation}\label{40}
\mathscr{K}^{(1)}_{u_0,h}(t)=\|((H_h-z)\phi(H_h)\mathcal{J}_{\pm,h}\phi(H_0)-\phi(H_h) \mathcal{J}_{\pm,h}(H_0-z)\phi(H_0))U_0(t)u_0\|_Y,
\end{equation}
to prove the assertion of the lemma, it is enough to prove that
\begin{equation}\label{41}
\mathcal{K}_{u_0,h}(t)=\|(\phi(H_h)\mathcal{J}_{\pm,h}-\mathcal{J}_{\pm,h} \phi(H_0))U_0(t)u_0\|_Y
\end{equation}
belongs to $L^1((-\infty,\infty);dt)$. To this end, for all $u_0,u\in X$ we have
\begin{equation}\label{42}
\langle \mathcal{J}_{\pm,h}u_0, H_hu\rangle-\langle \mathcal{J}_{\pm,h}H_0u_0 ,u\rangle= \langle G_{0,h}u_0,G_hu\rangle.
\end{equation}
By (\ref{42}) we have for any $v_0,v\in Y$
\begin{eqnarray}\label{43}
\begin{array}{ll}
\langle G_{0,h}R_0(z)v_0,G_hR_h(\overline{z})v\rangle&\!\!\!=\langle \mathcal{J}_{\pm,h}R_0(z)v_0,H_hR_h(\overline{z})v\rangle-\langle \mathcal{J}_{\pm,h}H_0R_0(z)v_0,R_h(\overline{z})v\rangle\\
&\!\!\!=\langle \mathcal{J}_{\pm,h}R_0(z)v_0,v\rangle-\langle \mathcal{J}_{\pm,h}v_0,R_h(\overline{z})v\rangle+\\
&\!\!\!+\langle \mathcal{J}_{\pm,h}R_0(z)v_0,\overline{z} R_h(\overline{z})v\rangle-\langle z\mathcal{J}_{\pm,h}R_0(z)v_0,R_h(\overline{z})v\rangle.
\end{array}
\end{eqnarray}

\begin{Theo}
\emph{
Given self-adjoint operators $T_h$ and $T_0$ in Hilbert spaces $\mathscr{H}$ and $\mathscr{H}_0$ respectively, let the WOs $W_{\pm}(T_h,T_0;J_h)$ and $W_{\pm}(T_0,T_h;J_h^*)$ exist, where $J_h$ is some bounded identification. Assume that, as $h\to\infty$, the operator $T_h$ is convergent in the SRS to $T_\infty$ and that $J_h$ and $J_h^*$ converge strongly to $J_\infty$ and $J_\infty^*$ respectively. If, for $T_h$, $T_0$, $\mathscr{H}$, $\mathscr{H}_0$ and $J_h$, the functions $\mathscr{K}^{(1)}_{u_0,h}(t)$ and $\mathscr{K}^{(2)}_{u,h}(t)$ satisfy the conclusions of Corollaries 1 and 2 respectively, then
\begin{equation}\label{47}
s\!-\!\!\!\dis\lim_{h\to\infty}W_{\pm}(T_h,T_0;J_h)= W_{\pm}(T_\infty ,T_0;J_\infty)
\end{equation}
and
\begin{equation}\label{48}
s\!-\!\!\!\dis\lim_{h\to\infty}W_{\pm}(T_0,T_h;J_h^*)=W_{\pm}(T_0 ,T_\infty;J_\infty^*).
\end{equation}
}
\end{Theo}
\hspace{-4mm}\underline{\emph{Proof}}. See Theorem 2.1 and Remark 2.3 in \cite{BRU}. \hfill{$\blacksquare$}\\

Since, in general, the strong convergence of an operator does not imply the strong convergence of its adjoint to the adjoint of its strong limit, therefore we have assumed the strong convergence of $J_h^*$ to $J_\infty^*$ parallel to the strong convergence of $J_h$ to $J_\infty$ in Theorem 3. Fortunately, the identification operator $\mathcal{J}_{\pm,h}$  is a PSDO with adjoint operator $\mathcal{J}^*_{\pm,h}$ given by
\begin{equation}\label{49}
(\mathcal{J}^*_{\pm,h}g)(\zeta)=(2\pi)^{-3/2}\dis\int_{\mathbb{R}^3}\dis e^{-i x\cdot\zeta - i\Phi_{\pm,h}(x,\zeta)} \mathcal{P}_{\pm,h}(x,\zeta)\mathscr{C}_{\pm}(x,\zeta)\psi(|\zeta|^2)g(x) \,dx.
\end{equation}
This implies that if $\mathcal{J}_{\pm,h}$ is strongly convergent to $\mathcal{J}_{\pm,\infty}$ as $h\to\infty$, then $\mathcal{J}^*_{\pm,h}$ is also strongly convergent to $\mathcal{J}^*_{\pm,\infty}$.\\

\subsubsection{The case $\rho\in(1/2,1)$}
According to the definition of the phase function $\Phi_{\pm,h}(x,\zeta)$, if $\rho\in(1/2,1)$, then $N=1$ satisfies the condition $(N+1)\rho>1$. In this case $\Phi_{\pm,h}(x,\zeta)$ can be chosen as, after neglecting the quadratic terms,
\begin{equation}\label{65}
\Phi_{\pm,h}(x,\zeta)= \pm\eta(\zeta)\dis\int_0^\infty(V_h(x\pm t\zeta)-V_h(\pm t\zeta))\,dt.
\end{equation}
By Proposition 3, to study the strong convergence of $W_{\pm,h}$ with the identification $\mathcal{J}_{\pm,h}$ (respectively $W_{\pm,h}^*$ with $\mathcal{J}_{\pm,h}^*$) is equivalent to study its strong convergence with $\boldsymbol{\mathcal{J}}^{(3)}_{\pm,h}$ (respectively with $\boldsymbol{\mathcal{J}}^{(3),*}_{\pm,h}$, where $\boldsymbol{\mathcal{J}}^{(3),*}_{\pm,h}$ is the adjoint operator of $\boldsymbol{\mathcal{J}}^{(3)}_{\pm,h}$),
\begin{equation}\label{66}
(\boldsymbol{\mathcal{J}}^{(3)}_{\pm,h}g)(x)=(2\pi)^{-3/2}\dis\int_{\mathbb{R}^3}\dis e^{i x\cdot\zeta + i\Phi_{\pm,h}(x,\zeta)} p_{0}(\zeta)\mathscr{C}_{\pm}(x,\zeta)\psi(|\zeta|^2)\hat{g}(\zeta) \,d\zeta,
\end{equation}
where $\Phi_{\pm,h}$ is given by (\ref{65}). By dominated convergence theorem, and since the integrand in (\ref{66}) is bounded, then if $\Phi_{\pm,h}(x,\zeta)$, given by (\ref{65}), converges to $\Phi_{\pm,\infty}(x,\zeta)$ in the SRS, then the identification $\boldsymbol{\mathcal{J}}^{(3)}_{\pm,h}$ is strongly convergent to $\boldsymbol{\mathcal{J}}^{(3)}_{\pm,\infty}$, where
\begin{equation}\label{67}
(\boldsymbol{\mathcal{J}}^{(3)}_{\pm,\infty}g)(x)=(2\pi)^{-3/2}\dis\int_{\mathbb{R}^3}\dis e^{i x\cdot\zeta + i\Phi_{\pm,\infty}(x,\zeta)} p_{0}(\zeta)\mathscr{C}_{\pm}(x,\zeta)\psi(|\zeta|^2)\hat{g}(\zeta) \,d\zeta.
\end{equation}

\end{document}